\documentclass[preprint,12pt]{elsarticle}

\usepackage{amsmath}
\usepackage{cases}
\usepackage{color}
\usepackage{hyperref}
\usepackage{amssymb}
\usepackage{graphicx}
\usepackage{graphics}
\usepackage{subfigure}
\usepackage{mathrsfs}
\usepackage{textcomp}
\usepackage{bm}
\usepackage{url}
\usepackage{CJK}
\usepackage{mathrsfs}

\makeatletter
  \newcommand\figcaption{\def\@captype{figure}\caption}
  \newcommand\tabcaption{\def\@captype{table}\caption}
\makeatother

\biboptions{numbers,sort&compress}

\begin{document}
\begin{frontmatter}



\title{On the Generalized Wavelet-Galerkin Method}


\author[mymainaddress]{Zhaochen Yang}
\author[mymainaddress,my3thaddress]{Shijun Liao\corref{cor1}}
\ead{sjliao@sjtu.edu.cn} \cortext[cor1]{Corresponding author.}

\address[mymainaddress]{School of Naval Architecture, Ocean and Civil Engineering, Shanghai Jiao Tong University, Shanghai 200240, China}


\address[my3thaddress]{Ministry-of-Education Key Laboratory in Scientific and Engineering Computing, Shanghai 200240, China}


\begin{abstract}

In the frame of the traditional wavelet-Galerkin method based on the compactly supported wavelets, it is important to calculate the so-called connection coefficients that are some integrals whose integrands involve products of wavelets, their derivatives as well as some known coefficients in considered differential equations.  However,  even for linear differential equations with non-constant coefficient,  the  computation of connect coefficients becomes rather time-consuming and often even impossible.  In this paper, we propose a generalized wavelet-Galerkin method based on the compactly supported wavelets, which is computationally very efficient even for differential equations with non-constant coefficients, {     no matter linear or nonlinear problems}.   Some related mathematical theorems are proved, based on which the basic ideas of the generalized wavelet-Galerkin method are described in details.  In addition, some examples are used to illustrate its validity and high efficiency. {    A nonlinear example shows that the generalized wavelet-Galerkin method is not only valid to solve nonlinear problems, but also possesses the ability to find new solutions of multi-solution problems.} This method can be widely applied to various types of {    both linear and nonlinear} differential equations in science and engineering.
\end{abstract}


\end{frontmatter}


\section{Introduction}

Wavelet, as a powerful mathematical tool, has been widely applied to solve many problems in science and engineering. Since the breakthrough in 1988 when Daubechies \cite{daubechies1988orthonormal, daubechies1992ten} constructed an orthogonal, compactly supported wavelet, there has been an increasing interest in wavelet-based methods for ordinary and partial differential equations. Briefly speaking, there exist three kinds of wavelet-based methods: the wavelet finite element method \cite{WFEMko1995class, WFEMma2003study, WFEMhan2006spline}, the wavelet collocation method \cite{WCMbertoluzza1996wavelet, WCMvasilyev1996dynamically, WCMvasilyev2005adaptive} and the wavelet-Galerkin method \cite{WGMamaratunga1994wavelet, WGMavudainayagam2000wavelet, WGMwang2012wavelet, WGMliu2013wavelet, WGMzhang2015wavelet}, whereas the wavelet-Galerkin method has gained the widest attention due to its good convergence and stability characteristics \cite{gomes1996convergence, reginska2000stability, liu2010daubechies}.

To obtain convergent, accurate and stable solutions efficiently by the wavelet-Galerkin method, it is of great importance to choose a proper type of wavelet. Generally, there are some main considerations for the choice of wavelet: orthogonality, compact support, interpolation properties as well as high algebraic accuracy and the ability to represent functions at different levels of resolution \cite{chen1996computation, liu2014Phd}. However, it is well known that compactly supported wavelets, such as Daubechies' wavelet \cite{daubechies1988orthonormal, daubechies1992ten}, Coiflet wavelet \cite{beylkin1991fast, huang2002coiflet} and the generalized Coiflet-type wavelet \cite{wang2001generalized, liu2014Phd}, have no analytical expressions. More seriously, the derivatives of these compactly supported wavelets are highly oscillatory.

Therefore, great contradictions and challenges arise. On the one hand, the wavelet-Galerkin method possesses great advantages in solving differential equations so that these compactly supported wavelets are expected to be used. On the other hand, due to the lack of analytical expressions of these compactly supported wavelets and high oscillation of their derivatives, there has been a great challenge to calculate the connection coefficients \cite{chen1996computation, reginska2000stability, bulut2016alternative}, which play an important role in the wavelet-Galerkin method for differential equations. Simply speaking, the connection coefficients are some integrals whose integrands involve products of the wavelet bases, their derivatives as well as other known functions. However, as mentioned above, for compactly supported wavelets, it is difficult and sometimes even unable to calculate these connection coefficients by the traditional numerical methods.

To tackle the problem, a lot of works have been done. Motivated by the successful construction of orthogonal interval wavelets \cite{chui1992wavelets, cohen1993wavelets}, Chen et al. \cite{chen1996computation} developed a method to compute a series of connection coefficients. In 2007, Zhang et al. \cite{zhang2007comments} corrected some errors of this method,  and develop better procedures for the computation of the connection coefficients. Similar to Chen et al. \cite{chen1996computation}, Wang \cite{wang2001generalized} developed a high-precision approach to calculate the connection coefficients for the generalized Coiflet-type wavelets.

However, the problem is not completely solved, because these works only provide some algorithms to calculate connection coefficients for special cases. If a differential equation contains some non-constant coefficients, some new problems arise. For example, let us consider a differential equation on the bounded interval $[0,1]$, which contains a term $\sin (x) \frac{d^nu(x)}{dx^n}$, where $\sin(x)$ is a non-constant coefficient and $u(x)$ is the unknown function to be solved. Suppose that the wavelet base functions are constructed as $\left\{\varphi_k(x)\right\}_{k=0}^{N}$, and an arbitrary function $u(x) \in L^2[0,1]$ can be approximated by
\begin{equation} \label{General eg:u(x)}
  u(x)\approx \sum_{k=0}^N \mu_k \varphi_k(x),
\end{equation}
where $\mu_k$, $k=0,1,2,...,N$ are unknown coefficients to be determined. To solve the differential equation, we substitute this approximation into the original equation. As a result, it contains the term
\begin{equation}
  \sin(x)\frac{d^nu(x)}{dx^n} \approx \sum_{k=0}^N \mu _k \sin(x) \frac{d^n\varphi_k(x)}{dx^n}.
\end{equation}
 Therefore, when the wavelet-Galerkin method is adopted, it is necessary to calculate the following connection coefficients
\begin{equation} \label{General eg:Gamma_kl}
  \int_0^1 \sin(x) \frac{d^n\varphi_k(x)}{dx^n} \varphi_l(x) dx,
\end{equation}
where $k,l=0,1,2,...,N$.
Obviously, it is extremely difficult and even impossible to calculate such kinds of integrals for compactly supported wavelets by known methods yet, especially when $\sin(x)$ is replaced by an arbitrary function $f(x)$. Note that such kind of connection coefficients are dependent upon the function $f(x)$ and thus one had to calculate them for different equations.   In theory, it is convenient if all connection coefficients are only dependent upon the wavelet base functions, i.e. independent of differential equations under consideration.

Of course, we can further approximate the coefficient function $\sin(x)$ by
\begin{equation} \label{General eg:sin(x)}
  \sin(x)\approx \sum_{k=0}^N \alpha_k \varphi_k(x),
\end{equation}
where $\alpha_k=\int_0^1 \sin(x) \varphi_k(x) dx$, $k=0,1,2,...,N$, so that the term $\sin(x) \frac{d^n u(x)}{dx^n}$ becomes
\begin{equation}
  \sin(x)\frac{d^n u(x)}{dx^n} \approx \sum_{k=0}^N \sum_{l=0}^N \mu_k \alpha_l \frac{d^n\varphi_k(x)}{dx^n} \varphi_l(x).
\end{equation}
In this way,  the connection coefficients are indeed independent of differential  equations under consideration.   However,  it is time-consuming to calculate the connection coefficients
\begin{equation} \label{General eg:Gamma_klm}
  \Gamma_{k,l,m}^n=\int_0^1 \frac{d^n\varphi_k(x)}{dx^n} \varphi_l(x) \varphi_m(x) dx,
\end{equation}
where $k,l,m=0,1,2,...,N$. In addition, it may also cause serious increase of the number of necessary terms, which is another reason that makes the algorithm seriously inefficient. Thus, connection coefficients expressed in the forms of (\ref{General eg:Gamma_kl}) and (\ref{General eg:Gamma_klm}) are called the complicated connection coefficients in this paper.

The problem of the above-mentioned  complicated connection coefficients seriously limits the wide applications of the wavelet-Galerkin method for differential equations with non-constant coefficients. Therefore, we are motivated to develop a generalized wavelet-Galerkin method which possesses the following advantages:
\begin{enumerate}
  \item [1)] The above-mentioned complicated connection coefficients is completely avoided;
  \item [2)] All connection coefficients are independent of differential equations under consideration, so that the algorithm is very efficient and adaptive for different boundary value problems no matter what kinds of non-constant coefficients are contained;
  {    
  \item [3)] It is valid not only for linear differential equations with non-constant coefficients, but also for nonlinear ones;
  }
  \item [4)] It can be conveniently extended to high-dimensional boundary value problems.
\end{enumerate}

In this paper, a generalized approach is developed to efficiently handle {    both linear and nonlinear} differential equations with non-constant coefficients in the frame of the wavelet-Galerkin method. The paper is arranged as follows: In Section 2, some properties of the generalized Coiflet-type wavelet are introduced. In Section 3, a few of  related mathematical theorems are given and proved. In Section 4, some linear cases with exact solutions are given to illustrate the validity, accuracy and computational efficiency of the generalized wavelet-Galerkin method. {     In Section 5, a comparatively complicated nonlinear application is presented to show the validity as well as the potential of the generalized wavelet-Galerkin method. In Section 6,} some concluding remarks are given.

\section{Some properties of the generalized Coiflet-type wavelet}

The generalized Coiflet-type wavelet is defined by its scaling function $\varphi(x)$ and wavelet function $\psi(x)$ which possess the following properties \cite{wang2001generalized, liu2014wavelet}:
\begin{eqnarray} \label{Eq:WaveProp}
    &&(a)\ \varphi(x)=\sum_{k \in \mathbb{Z}}p_{k} \varphi (2x-k);  \\
    &&(b)\ \psi(x)=\sum_{k \in \mathbb{Z}}(-1)^k p_{1-k} \psi (2x-k);\\
    &&(c)\ M_n=M_1^n,\quad for\  0\leq n<N;\\
    &&(d)\ \int_{-\infty}^{+\infty}x^n \psi(x)dx=0, \quad for \ 0\leq n<N; \qquad \qquad \qquad \qquad \qquad \qquad \\
    &&(e)\ \sum_{k \in \mathbb{Z}}p_{k} \varphi (x-k)=1;
\end{eqnarray}
where $p_k$  is the low-pass filter coefficients \cite{liu2014Phd, zhang2013influence}, $N$ is the number of vanishing moment, $M_n=\int_{-\infty}^{+\infty}x^n \varphi(x)dx$ is the nth-order moment of the scaling function and the corresponding supported interval is $supp\left[\varphi(x)\right]=\left[0, \ 3N-1\right]$. Without loss of generality, the generalized Coiflet-type wavelet with $N=6, M_1=7$ is used in this paper.

The generalized Coiflet-type wavelet possesses a lot of advantages: orthogonality, compact support, quasi-interpolation (not exactly interpolation but with high accuracy) properties as well as high algebraic accuracy and the ability to represent functions at different levels of resolution. Based on these advantages, wavelet basis for $L^2$-function space is constructed by Wang et al. \cite{wang2001generalized, zhou2011wavelet}, and an arbitrary function $f(x) \in L^2[0,1]$ can be approximated by
\begin{equation}\label{Def:WaveletAppr:1D}
    f(x) \approx P^j_x f(x) = \sum_{k=0}^{2^j} f\left(\frac{k}{2^j}\right) \varphi_{j,k}(x),
\end{equation}
in which $P^j_x$ is an projective operator with respect to $x$, $j$ is the resolution level and $\varphi_{j,k}(x)$ is the wavelet basis defined by
\begin{equation} \label{Def:WaveletBasis}
\varphi_{j,k}(x)=\left\{
\begin{aligned}
         &\sum_{i=2-3N+M_1}^{-1}T_{0,k}\left(\frac{i}{2^j}\right)\varphi(2^j x- i+M_1)+\varphi(2^jx-k+M_1), \\
         &\qquad \qquad \qquad \qquad \qquad \qquad \qquad \qquad \qquad k\in [0,3], \\
         &\varphi(2^jx-k+M_1), \\
         &\qquad \qquad \qquad \qquad \qquad \qquad \qquad \qquad \qquad k\in [4,2^j-4], \\
         &\sum_{i=2^j+1}^{2^j-1+M_1}T_{1,2^j-k}\left(\frac{i}{2^j}\right)\varphi(2^j x- i+M_1)+\varphi(2^jx-k+M_1), \\
         &\qquad \qquad \qquad \qquad \qquad \qquad \qquad \qquad \qquad k\in [2^j-3,2^j],\\
\end{aligned} \right.
\end{equation}
where
\begin{equation}\label{Def:FunT_il.}
    T_{0,k}(x)=\sum_{i=0}^3 \frac{p_{0,i,k}}{i!}x^i \ , \ T_{1,k}(x)=\sum_{i=0}^3 \frac{p_{1,i,k}}{i!}(x-1)^i,
\end{equation}
and coefficients $p_{0,i,k}$ and $p_{1,i,k}$ are assigned as
\begin{eqnarray}\label{Def:Coef.P}
    \mathbf{P_0}&=&\left[
  \begin{array}{cccc}
    1 & 0 & 0 & 0\\
    -11/6 &  3 & -3/2 & 1/3\\
    2 & -5 &  4 & -1\\
    -1 & 3 & -3 & 1\\
  \end{array}
  \right],  \\
   \mathbf{P_1}&=&\left [
  \begin{array}{cccc}
    1 & 0 & 0 & 0\\
    11/6 &  -3 & 3/2 & -1/3\\
    2 & -5 &  4 & -1\\
    1 & -3 & 3 & -1\\
  \end{array}
  \right ]
\end{eqnarray}
through the relations $\mathbf{P_0}=\{2^{-ij}p_{0,i,k}\}$ and $\mathbf{P_1}=\{2^{-ij}p_{1,i,k}\}$ with $i,k=0,1,2,3$ and $j$ being the resolution level.

Actually, $p_{0,i,k}$ and $p_{1,i,k}$ represent the difference coefficients of $f^{(i)}(0)$ and $f^{(i)}(1)$, respectively. Let $h=\frac{1}{2^j}$, then according to Taylor's theorem, $f(kh)$  and $f(1-kh)$ can be expanded into
\begin{equation}\label{Eq:f(kh)}
    f(kh)=f(0)+f'(0)(kh)+\frac{f''(0)}{2!}(kh)^2+\frac{f''(0)}{3!}(kh)^3+o(h^3),
\end{equation}
and
\begin{equation}\label{Eq:f(1-kh)}
    f(1-kh)=f(1)+f'(1)(-kh)+\frac{f''(1)}{2!}(-kh)^2+\frac{f''(1)}{3!}(-kh)^3+o(h^3),
\end{equation}
respectively. Therefore, the coefficients $p_{0,i,k}$ and $p_{1,i,k}$  can be obtained by the following equations
\begin{equation}\label{Eq:d=Pf}
    \mathbf{d_0}=\mathbf{P_0}\mathbf{f_0}, \ \ \mathbf{d_1}=\mathbf{P_1}\mathbf{f_1},
\end{equation}
where $\mathbf{d_0}=\{f(0),hf'(0),h^2f''(0),h^3f'''(0)\}^T$, $\mathbf{f_0}=\{f(0),f(h),f(2h),f(3h)\}^T$, $\mathbf{d_1}=\{f(1),hf'(1),h^2f''(1),h^3f'''(1)\}^T$, and $\mathbf{f_1}=\{f(1),f(1-h),f(1-2h),f(1-3h)\}^T$,  $f(kh)$  and $f(1-kh)$  are approximated by Eqs.(\ref{Eq:f(kh)}) and (\ref{Eq:f(1-kh)}) without the high-order term $o(h^3)$.

Obviously, the derivative of the function can be approximated by
\begin{equation}\label{Def:WaveletAppr:1D:Dn}
    \frac{d^nf(x)}{dx^n} \approx \frac{d^nP^j_x f(x)}{dx^n} = \sum_{k=0}^{2^j} f\left(\frac{k}{2^j}\right) \frac{d^n\varphi_{j,k}(x)}{dx^n},
\end{equation}
where the values of $\frac{d^n\varphi_{j,k}(x)}{dx^n}$ on equinoxes can be calculated through the method suggested by Wang \cite{wang2001generalized}.

Of course, for boundary value problems, there are always some boundary conditions to be satisfied. So some minor modifications need to be made on the approximation (\ref{Def:WaveletAppr:1D}) according to different boundary conditions.  For example,
\begin{enumerate}
  \item [(a)] When the value of $f(0)$ or $f(1)$ is given, the boundary value condition can be satisfied by directly substituting these two given values into wavelet approximation (\ref{Def:WaveletAppr:1D}).
  \item [(b)] When $f'(0)=0$, the wavelet basis should be modified into
    \begin{equation}\label{Def:wavelet modification:tilde}
        \tilde{\varphi}_{j,k}(x)=\varphi_{j,k}(x)|_{p_{0,1,k}\rightarrow0}.
    \end{equation}
  \item [(c)] When $f'(0)=0$ and $f''(0)=0$, the wavelet basis should be modified into
    \begin{equation}\label{Def:wavelet modification:double-tilde}
     \tilde{ \tilde{\varphi}}_{j,k}(x)=\varphi_{j,k}(x)|_{p_{0,1,k}\rightarrow0, \ p_{0,2,k}\rightarrow0}.
    \end{equation}
\end{enumerate}
For more details on the modifications, please refer to \cite{liu2014Phd}.

According to \cite{liu2014Phd, sweldens1994quadrature}, the error estimation of the wavelet approximation (\ref{Def:WaveletAppr:1D}) is
{    
\begin{equation}\label{Eq:AccuracyPj:1D}
    \left\| \frac{d^nf(x)}{dx^n}-\frac{d^nP^j_x f(x)}{dx^n}\right\| _{L^2} \leq C 2^{-j(N-n)},
\end{equation}
}
where $C$ is a positive constant that depends only on the function $f(x)$ and low-pass filter coefficients $p_k$, $N$ is the number of vanishing moment, and $0 \leq n <N$.

According to the theory of multiresolution analysis \cite{meyer1995wavelets}, two-dimensional wavelet basis can be directly extended by the tensor products of one-dimensional wavelet basis. Therefore, for an arbitrary function $f (x, y) \in L^2[0, 1]^2$, we have the following approximations
\begin{equation}\label{Def:WaveletAppr:2D}
    f(x,y) \approx P^j_x P^j_y f(x,y) = \sum_{k=0}^{2^j}\sum_{l=0}^{2^j} f\left(\frac{k}{2^j}, \frac{l}{2^j}\right) \varphi_{j,k}(x)\varphi_{j,l}(y)
\end{equation}
and
{    
\begin{equation}\label{Def:WaveletAppr:1D:Dn}
    \frac{\partial^{m+n}f(x,y)}{\partial x^m \partial y^n} \approx \frac{\partial^{m+n}P^j_x P^j_y f(x,y)}{\partial x^m \partial y^n} = \sum_{k=0}^{2^j}\sum_{l=0}^{2^j} f\left(\frac{k}{2^j}, \frac{l}{2^j}\right) \frac{d^m \varphi_{j,k}(x)}{dx^m}\frac{d^n \varphi_{j,l}(y)}{dy^n}.
\end{equation}
}
The error estimation can be derived from Eq. (\ref{Eq:AccuracyPj:1D}) as
{    
\begin{equation}\label{Eq:AccuracyPj:2D}
    \left\| \frac{\partial^{m+n}f(x,y)}{\partial x^m \partial y^n}-\frac{\partial^{m+n}P^j_x P^j_yf(x,y)}{\partial x^m \partial y^n}\right\| _{L^2} \leq C 2^{-\min \{j(N-m), j(N-n)\}}.
\end{equation}
}

\section{Some related theorems}

In this section, some fundamental theorems are proposed and proved. Based on these theorems, the generalized wavelet-Galerkin method for linear boundary value problems governed by differential equations with non-constant coefficients is deduced. {     For simplicity, we mainly describe the basic ideas as well as the theorems with linear cases, while all of them are similar for nonlinear ones.
}

\subsection{Some basic theorems and definitions}
\textbf{Theorem 1} Assuming that $I$ is a bounded open interval, $n$ is a positive integer, $f(x) \in C^{n}(I)$ and $u(x) \in C^{n}(I)$, then in the interval $I$ it holds
\begin{equation}\label{Th1:mainEq}
     f(x)\frac{d^n u(x)}{dx^n}=\sum_{k=0}^{n}(-1)^k\binom{n}{k}\frac{d^{n-k}}{dx^{n-k}}\left[\frac{d^kf(x)}{dx^k}u(x)\right],
\end{equation}
where the $0^{th}$-order derivative of arbitrary function $g(x) \in C(I)$ is defined as itself, i.e.
\begin{equation}\label{Th1:def:0order}
   \frac{d^0 g(x)}{d x^0 }\triangleq g(x).
\end{equation}

\textbf{\emph{Proof}}  Considering $n \in \mathbb{Z}^+$, it is equivalent to dealing with the following open statement:
\begin{equation}
    S(n):  f(x)\frac{d^n u(x)}{dx^n}=\sum_{k=0}^{n}(-1)^k\binom{n}{k}\frac{d^{n-k}}{dx^{n-k}}\left[\frac{d^kf(x)}{dx^k}u(x)\right], \nonumber
\end{equation}
for arbitrary function $f(x), u(x) \in C^n(I)$.

Here, we prove this theorem by the method of mathematical induction.

\textit{Basic Step}: we start with the statement $S(1)$ and find that
\begin{align}
    S(1):  f(x)\frac{d u(x)}{dx}&=\frac{d}{dx}\left[f(x)u(x)\right]-\frac{df(x)}{dx}u(x), \nonumber \\
    &=\sum_{k=0}^{1}(-1)^k\binom{1}{k}\frac{d^{1-k}}{dx^{1-k}}\left[\frac{d^kf(x)}{dx^k}u(x)\right], \nonumber
\end{align}
holds for arbitrary $f(x), u(x) \in C^1(I)$, say, $S(1)$ is true.

\textit{Inductive Step}: Assuming the truth of $S(k)$ for some particular $k \in \mathbb{Z}^+$, we have
\begin{equation}
    S(k):  f(x)\frac{d^k u(x)}{dx^k}=\sum_{i=0}^{k}(-1)^i\binom{k}{i}\frac{d^{k-i}}{dx^{k-i}}\left[\frac{d^if(x)}{dx^i}u(x)\right], \nonumber
\end{equation}
for arbitrary $f(x), u(x) \in C^k(I)$.

Then, for arbitrary $f(x), u(x) \in C^{k+1}(I)$, it is true that $\frac{df(x)}{dx} \in C^{k}(I)$ and $\frac{d^kf(x)}{dx^k} \in C^{1}(I)$. Therefore,
\begin{align}
    f(x)\frac{d^{k+1} u(x)}{dx^{k+1}}&=f(x)\frac{d}{dx}\left[\frac{d^{k} u(x)}{dx^{k}}\right] \nonumber \\
    &=\frac{d}{dx}\left[f(x)\frac{d^{k} u(x)}{dx^{k}}\right] -\frac{df(x)}{dx}\frac{d^{k} u(x)}{dx^{k}} \nonumber \\
    &=\frac{d}{dx} \left\{\sum_{i=0}^{k}(-1)^i\binom{k}{i}\frac{d^{k-i}}{dx^{k-i}}\left[\frac{d^if(x)}{dx^i}u(x)\right]\right\} \nonumber\\
    &\quad-\sum_{i=0}^{k}(-1)^i\binom{k}{i}\frac{d^{k-i}}{dx^{k-i}}\left[\frac{d^i}{dx^i}\left( \frac{df(x)}{dx} \right) u(x)\right].  \nonumber
 \end{align}
Using the induction hypothesis, we have
\begin{align}
    f(x)\frac{d^{k+1} u(x)}{dx^{k+1}}&=\sum_{i=0}^{k}(-1)^i\binom{k}{i}\frac{d^{k+1-i}}{dx^{k+1-i}}\left[\frac{d^if(x)}{dx^i}u(x)\right] \nonumber \\
    &\quad-\sum_{i=0}^{k}(-1)^i\binom{k}{i}\frac{d^{k-i}}{dx^{k-i}}\left[\frac{d^{i+1}f(x)}{dx^{i+1}}\ u(x)\right]  \nonumber\\
    &=\sum_{i=0}^{k}(-1)^i\binom{k}{i}\frac{d^{k+1-i}}{dx^{k+1-i}}\left[\frac{d^if(x)}{dx^i}u(x)\right] \nonumber \\
    &\quad+\sum_{j=1}^{k+1}(-1)^{j}\binom{k}{j-1}\frac{d^{k+1-j}}{dx^{k+1-j}}\left[\frac{d^{j}f(x)}{dx^{j}}\ u(x)\right]  \nonumber \\
    &=(-1)^0\binom{k}{0}\frac{d^{k+1}}{dx^{k+1}}\left[\frac{d^0f(x)}{dx^0}u(x)\right]+(-1)^{k+1}\binom{k}{k}\frac{d^{0}}{dx^{0}}\left[\frac{d^{k+1}f(x)}{dx^{k+1}}\ u(x)\right] \nonumber \\
    &\quad+\sum_{i=1}^{k}(-1)^i \left[\binom{k}{i}+\binom{k}{i-1}\right] \frac{d^{k+1-i}}{dx^{k+1-i}}\left[\frac{d^if(x)}{dx^i}u(x)\right].\nonumber
 \end{align}

Since $\binom{k+1}{0}=\binom{k}{0}=\binom{k+1}{k+1}=\binom{k}{k}=1$ and $\binom{k}{i}+\binom{k}{i-1}=\binom{k+1}{i}$, we have
\begin{equation}\label{Pf:Sk_1_2}
    f(x)\frac{d^{k+1} u(x)}{dx^{k+1}}=\sum_{i=0}^{k+1}(-1)^i \binom{k+1}{i} \frac{d^{k+1-i}}{dx^{k+1-i}}\left[\frac{d^if(x)}{dx^i}u(x)\right],\nonumber
\end{equation}
which indicates that the statement $S(k+1)$ is true.

Hence, for any particular  $k \in \mathbb{Z}^+$, it follows that $S(k) \Longrightarrow S(k+1)$. Considering the truth of $S(1)$, it leads to the conclusion that the statement $S(n)$ is true  for all $n \in \mathbb{Z}^+$.

This ends the proof. \\

\textbf{Theorem 2} Let $I$ denote a bounded open interval. The linear ordinary differential equation
\begin{equation}\label{Th2:ODE}
    \sum_{n=0}^N a_{n}(x)\frac{d^{n}u(x)}{d x^n}=r(x), \ \ x \in I
\end{equation}
where $a_n(x) \in C^{n}(I)$ and $r(x)$ are known functions can be rewritten in the form
\begin{equation}\label{Th2:ODE:trans}
   \sum_{n=0}^N  \frac{d^{n}\left(b_{n}(x)u(x)\right)}{d x^n}=r(x), \ \ x \in I,
\end{equation}
where
\begin{equation}\label{Th2:ODE:coeff_fun}
   b_{n}(x)=\sum_{r=n}^N \left(-1\right)^{r-n}\binom{r}{n}\frac{d^{r-n}a_{r}(x)}{d x^{r-n}}.
\end{equation}

\textbf{\emph{Proof}} According to \textbf{Theorem 1}, for each $a_n(x) \in C^n(I)$, $n=0,1,2,...N$, we have
\begin{align}
a_n(x)\frac{d^n u(x)}{dx^n}&=\sum_{k=0}^{n}(-1)^k\binom{n}{k}\frac{d^{n-k}}{dx^{n-k}}\left[\frac{d^ka_n(x)}{dx^k}u(x)\right] \nonumber \\
&=\sum_{k=0}^{n}(-1)^{n-k}\binom{n}{k}\frac{d^{k}}{dx^{k}}\left[\frac{d^{n-k}a_n(x)}{dx^{n-k}}u(x)\right]. \nonumber
\end{align}

Therefore,
\begin{align}
    \sum_{n=0}^N a_{n}(x)\frac{d^{n}u(x)}{d x^n}&=\sum_{n=0}^N\sum_{k=0}^{n} (-1)^{n-k}\binom{n}{k}\frac{d^{k}}{dx^{k}}\left[\frac{d^{n-k}a_n(x)}{dx^{n-k}}u(x)\right] \nonumber \\
    &=\sum_{k=0}^N\sum_{n=k}^{N}(-1)^{n-k}\binom{n}{k}\frac{d^{k}}{dx^{k}}\left[\frac{d^{n-k}a_n(x)}{dx^{n-k}}u(x)\right]. \nonumber
\end{align}
Rewriting the summation index $n$ as $r$, and $k$ as $n$ on the right side, we have
\begin{align}
    \sum_{n=0}^N a_{n}(x)\frac{d^{n}u(x)}{d x^n}&=\sum_{n=0}^N\sum_{r=n}^{N}(-1)^{r-n}\binom{r}{n}\frac{d^{n}}{dx^{n}}\left[\frac{d^{r-n}a_r(x)}{dx^{r-n}}u(x)\right]  \nonumber \\
    &=\sum_{n=0}^N \frac{d^{n}}{dx^{n}} \left\{ \left[\sum_{r=n}^{N}(-1)^{r-n}\binom{r}{n}\frac{d^{r-n}a_r(x)}{dx^{r-n}}\right] u(x)\right\}. \nonumber
\end{align}\
This ends the proof.\\

\textbf{Theorem 3} Assuming that $\Omega$ is a bounded open region, $f(x,y) \in C^{m+n}(\Omega)$ and $u(x,y) \in C^{m+n}(\Omega)$, then in the region $\Omega$ it holds
\begin{equation}\label{Th3:mainEq}
     f(x,y)\frac{\partial^{m+n}u(x,y)}{\partial x^m \partial y^n}=\sum_{i=0}^{m}\sum_{j=0}^{n}(-1)^{i+j}\binom{m}{i}\binom{n}{j}\frac{\partial^{(m-i)+(n-j)}}{\partial x^{m-i} \partial y^{n-j}}\left[\frac{\partial^{i+j}f(x,y)}{\partial x^i \partial y^j}u(x,y)\right],
\end{equation}
where the derivation operators are defined by
\begin{equation}\label{Th3:def:0order}
    \left\{
    \begin{aligned}
    &\frac{\partial^{m+0}g(x,y)}{\partial x^m \partial y^0}\triangleq\frac{\partial^m g(x,y)}{\partial x^m} \\
    &\frac{\partial^{0+n}g(x,y)}{\partial x^0 \partial y^n}\triangleq\frac{\partial^n g(x,y)}{\partial y^n} \\
    &\frac{\partial^{0+0}g(x,y)}{\partial x^0 \partial y^0}\triangleq g(x,y).
    \end{aligned}
    \right.
\end{equation}

\textbf{Theorem 4} Let $\Omega$ denote a bounded open region. The linear partial differential equation
\begin{equation}\label{Th4:PDE}
    \sum_{m=0}^M \sum_{n=0}^N a_{m,n}(x,y)\frac{\partial^{m+n}u(x,y)}{\partial x^m \partial y^n}=r(x,y), \ \ (x,y) \in \Omega
\end{equation}
where $a_{m,n}(x,y) \in C^{m+n}(\Omega)$ and $r(x,y)$ are known functions can be rewritten in the form
\begin{equation}\label{Th4:PDE:trans}
   \sum_{m=0}^M \sum_{n=0}^N \frac{\partial^{m+n}\left(b_{m,n}(x,y)u(x,y)\right)}{\partial x^m \partial y^n}=r(x,y), \ \ (x,y) \in \Omega,
\end{equation}
where
\begin{equation}\label{Th4:PDE:coeff_fun}
   b_{m,n}(x,y)=\sum_{r=m}^M \sum_{s=n}^N \left(-1\right)^{(r-m)+(s-n)}\binom{r}{m}\binom{s}{n}\frac{\partial^{(r-m)+(s-n)}a_{r,s}(x,y)}{\partial x^{r-m} \partial y^{s-n}}.
\end{equation}

Then, let us define some operators that will be used later.\\

\textbf{Definition 1 (Hadamard/Schur product)} Let $\mathbf{A}=\left\{a_{kl}\right\}_{m\times n}$, $\mathbf{B}=\left\{b_{kl}\right\}_{m\times n}$ $\in \mathbb{R}^{m \times n}$ (vectors when $n=1$), then the Hadamard/Schur product of matrixes or vectors is defined by
\begin{equation}\label{Def:dot product}
   \mathbf{A\circ B}=\left\{c_{kl}=a_{kl}b_{kl}\right\}_{m\times n}.
\end{equation}

\textbf{Definition 2 (Column-dot product)} Let $\mathbf{A}=\left\{a_{kl}\right\}_{m\times n} \in \mathbb{R}^{ m \times n}$, $\mathbf{v}=\left\{v_{l}\right\}_{n \times 1} \in \mathbb{R}^{n\times 1}$ , then the column-dot product is defined by
\begin{equation}\label{Def:column-dot product}
   \mathbf{A}\odot \mathbf{v}=\left\{c_{kl}=a_{kl}v_{l}\right\}_{m\times n}.
\end{equation}

\subsection{The generalized wavelet-Galerkin method}
Using the above-mentioned theorems, we propose here the generalized wavelet-Galerkin method. Without loss of generality, let us consider the linear ordinary and partial differential equations (\ref{Th2:ODE}) and (\ref{Th4:PDE}), which are difficult to solve by means of the traditional wavelet-Galerkin method if $a_n(x)$ and $a_{m,n}(x,y)$ are non-constants.  However, according to the basic theorems, they can be rewritten as Eqs. (\ref{Th2:ODE:trans}) and (\ref{Th4:PDE:trans}), which are more convenient to solve by means of the wavelet-Galerkin method, as shown below. In this way, the connection coefficients in the frame of the wavelet-Galerkin method are independent of differential equations.  This is the basic idea of the generalized wavelet-Galerkin method, {    which is valid for both linear and nonlinear differential equations.}

\subsubsection{Approach for one-dimensional BVPs}
Let us consider the generalized linear boundary value problem (BVP) on the interval $[0,1]$ governed by the ordinary differential equation
\begin{equation}\label{Eq:ODE:trans}
   \sum_{n=0}^N  \frac{d^{n}\left(b_{n}(x)u(x)\right)}{d x^n}=r(x), \ \ 0<x<1
\end{equation}
with some boundary conditions, where $r(x)$ and $b_{n}(x)$ are known functions.

Instead of directly approximating the unknown function $u(x)$ by wavelet basis, in the frame of the generalized wavelet-Galerkin method we take $b_n(x) u(x)$ as a whole function. Using wavelet approximation (\ref{Def:WaveletAppr:1D}) and considering the boundary conditions, $b_n(x) u(x)$ can be approximated by
\begin{equation}\label{Eq:WaveletAppr:1D:fu}
    f_n(x)u(x) \approx P^j_x \left[f_n(x)u(x)\right] = \sum_{k=0}^{2^j} f_n\left(\frac{k}{2^j}\right) u\left(\frac{k}{2^j}\right) \varphi_{j,k}(x)
\end{equation}
for $n=1,2,...,N$,  and $r(x)$ can be approximated by
\begin{equation}\label{Eq:WaveletAppr:1D:r}
    r(x) \approx P^j_x \left[r(x)\right] = \sum_{k=0}^{2^j} r\left(\frac{k}{2^j}\right) \varphi_{j,k}(x).
\end{equation}

Substituting Eqs. (\ref{Eq:WaveletAppr:1D:fu}) and (\ref{Eq:WaveletAppr:1D:r}) into Eq. (\ref{Eq:ODE:trans}), we have
\begin{equation}\label{Eq:ODE:wavelet}
   \sum_{n=0}^N \sum_{k=0}^{2^j} b_n\left(\frac{k}{2^j}\right) u\left(\frac{k}{2^j}\right) \frac{d^n\varphi_{j,k}(x)}{dx^n}=\sum_{k=0}^{2^j} r\left(\frac{k}{2^j}\right) \varphi_{j,k}(x).
\end{equation}
Multiplying both sides of Eq. (\ref{Eq:ODE:wavelet}) by $\varphi_{j,l}(x)$ where $l=0, 1, 2, ..., 2^j$ (or from 1 to $2^j-1$ if values of $u(x)$ at $x=0$ and $x=1$ are given), respectively, and integrating over the interval $[0,1]$, it yields
\begin{equation}\label{Eq:ODE:wavelet-Galerkin:temp}
   \sum_{n=0}^N \mathbf{A_n}^T \left(\mathbf{b_n \circ u}\right)=\mathbf{A_0}^T\mathbf{r},
\end{equation}
where
\begin{equation}
    \left\{
    \begin{aligned}
    &\mathbf{A_n}=\left\{a_{kl}=\Gamma_{k,l}^{j,n}=\int_{0}^{1}\frac{d^n\varphi_{j,k}(\xi)}{d\xi^n} \varphi_{j,l}(\xi)d\xi\right\}_{(2^j+1)\times (2^j+1)}  \\
    & \mathbf{b_{n}}=\left\{b_{k}=b_{n}\left(\frac{k}{2^j}\right)\right\}_{(2^j+1)\times 1}\\
      & \mathbf{u}=\left\{u_{k}=u\left(\frac{k}{2^j}\right)\right\}_{(2^j+1)\times 1} \\
    & \mathbf{r}=\left\{r_{k}=r\left(\frac{k}{2^j} \right)\right\}_{(2^j+1)\times 1}.
    \end{aligned}
    \right.
\end{equation}
Eq. (\ref{Eq:ODE:wavelet-Galerkin:temp}) can be rewritten as
\begin{equation}\label{Eq:ODE:wavelet-Galerkin}
   \left[\sum_{n=0}^N \mathbf{A_n}^T \odot \mathbf{b_n}\right] \mathbf{u}=\mathbf{A_0}^T\mathbf{r}.
\end{equation}

In practice, there might be some minor differences due to the different boundary conditions, but the basic form of the wavelet-Galerkin equation does not change. Eq. (\ref{Eq:ODE:wavelet-Galerkin}) indicates that the connection coefficients are independent of $b_n(x)$ at all so that we only need to calculate the connection coefficients in the simplest form
\begin{equation}\label{Eq: Connection coefficients}
  \Gamma_{k,l}^{j,n}=\int_{0}^{1}\frac{d^n\varphi_{j,k}(\xi)}{d\xi^n} \varphi_{j,l}(\xi)d\xi,
\end{equation}
which is easy to gain in the frame of the wavelet-Galerkin method. Especially, since the connection coefficients are independent of differential equations, it is easy to build a database of the connection coefficients for data sharing.  This helps standardize the calculation procedure and can greatly increase the computational efficiency, as shown later by examples. Note that, according to \textbf{Theorem 2}, any a linear differential equation (\ref{Th2:ODE}) can be rewritten in the form of (\ref{Eq:ODE:trans}).  So, this approach is valid in general.

{    More importantly, in the frame of the generalized wavelet-Galerkin method, similar wavelet-Galerkin equations can be deduced for nonlinear problems, which is in the form of
\begin{equation}\label{Eq:NonlinearODE:wavelet-Galerkin}
   \left[\sum_{n=0}^N \mathbf{A_n}^T \odot \mathbf{b_n}\right] \mathbf{u}+ \left[\sum_{n=0}^N \mathbf{A_n}^T \odot \mathbf{c_n}\right] \mathbf{\mathscr{N}[u]}=\mathbf{A_0}^T\mathbf{r}.
\end{equation}
in which $\mathbf{c_n}$ and $\mathbf{\mathscr{N}[u]}$  are derived from the nonlinear part  with non-constant coefficients. An application will be given later.
}

\subsubsection{Approach for two-dimensional BVPs}
Then, let us consider the generalized linear BVP on the region $[0,1]^2$ governed by the partial differential equation
\begin{equation}\label{Eq:PDE:trans}
   \sum_{m=0}^M \sum_{n=0}^N \frac{\partial^{m+n}\left(b_{m,n}(x,y)u(x,y)\right)}{\partial x^m \partial y^n}=r(x,y), \ 0<x, y<1
\end{equation}
with some boundary conditions.

Using wavelet approximation (\ref{Def:WaveletAppr:2D}) and considering the boundary conditions, we similarly have the following equation
\begin{align}\label{Eq:PDE:wavelet}
   &\sum_{k=0}^{2^j}\sum_{l=0}^{2^j}\left[u\left(\frac{k}{2^j},\frac{k}{2^j}\right) \sum_{m=0}^M \sum_{n=0}^N b_{m,n}\left(\frac{k}{2^j},\frac{l}{2^j}\right)\frac{d^m\varphi_{j,k}(x)}{dx^m}\frac{d^n\varphi_{j,l}(y)}{dy^n}\right] \nonumber \\
   & \ \ \ \ \ \ \ \ \  \ \ \ \ \ \ \ \  \ \ \ \ \ \ \ \ \ \ \ \ \  \ \ \ \ =\sum_{k=0}^{2^j}\sum_{l=0}^{2^j}r\left(\frac{k}{2^j},\frac{l}{2^j}\right)\varphi_{j,k}(x)\varphi_{j,l}(y),
\end{align}
which leads to
\begin{equation}\label{Eq:PDE:wavelet-Galerkin:temp}
   \sum_{m=0}^M\sum_{n=0}^N \mathbf{A_m}^T \left(\mathbf{B_{mn}\circ U}\right) \mathbf{A_n}  =\mathbf{A_0}^T \mathbf{R} \mathbf{A_0}.
\end{equation}
where
\begin{equation} \label{Def:Matrixes:2D}
    \left\{
    \begin{aligned}
    &\mathbf{A_n}=\left\{a_{kl}=\Gamma_{k,l}^{j,n}=\int_{0}^{1} \frac{d^n\varphi_{j,k}(\xi)}{dx^n} \varphi_{j,l}(\xi)d\xi\right\}_{(2^j+1)\times (2^j+1)}  \\
    & \mathbf{B_{mn}}=\left\{b_{kl}=b_{m,n}\left(\frac{k}{2^j},\frac{l}{2^j}\right)\right\}_{(2^j+1)\times (2^j+1)}\\
      & \mathbf{U}=\left\{u_{kl}=u\left(\frac{k}{2^j},\frac{l}{2^j} \right)\right\}_{(2^j+1)\times (2^j+1)} \\
    & \mathbf{R}=\left\{r_{kl}=r\left(\frac{k}{2^j},\frac{l}{2^j} \right)\right\}_{(2^j+1)\times (2^j+1)}.
    \end{aligned}
    \right.
\end{equation}
Taking $rvec(\cdot)$ operator on both sides and using related operation properties, we obtain the following wavelet-Galerkin equation
 \begin{equation}\label{Eq:PDE:wavelet-Galerkin}
   \left[\sum_{m=0}^M\sum_{n=0}^N\left(\mathbf{A_m}\otimes \mathbf{A_n}\right)^T \odot rvec\left(\mathbf{B_{mn}}\right)\right] rvec\left(\mathbf{U}\right) =\left(\mathbf{A_0}\otimes \mathbf{A_0}\right)^T rvec\left(\mathbf{R}\right),
\end{equation}
{     where $\otimes$ denotes the Kronecker tensor product of matrixes and $rvec(\cdot)$ pulls a matrix into a vector row by row.}

 It can be found that the connection coefficients are independent of the coefficients $b_{m,n}(x,y)$ of the differential equation at all.  In other words,  the connection coefficients are independent of differential equations.  Note also that, according to \textbf{Theorem 4}, any two-dimensional differential equation (\ref{Th4:PDE}) can be rewritten as (\ref{Eq:PDE:trans}).  So, this approach is rather general.  {    Also, in the frame of the generalized wavelet-Galerkin method, similar wavelet-Galerkin equations can be deduced for nonlinear problems. In addition, both linear and nonlinear} boundary value problems in higher dimensions can also be solved in a similar way by means of the generalized wavelet-Galerkin method.

 {    
 In addition, as there is no truncation in \textbf{Theorem 1} $\sim$ \textbf{Theorem 4}, the error of the generalized wavelet-Galerkin method can be estimated by Eq. (\ref{Eq:AccuracyPj:1D}) and (\ref{Eq:AccuracyPj:2D}) only if $m, n$ are replaced by the order of the differential equations.
 }

\section{Illustrative linear examples}

To verify the validity and the computational efficiency of the generalized wavelet-Galerkin method, some examples of boundary value problems governed by linear ordinary and partial differential equations with different kinds of non-constant coefficients are solved. All results given in this section are obtained using a desktop (DELL Inspiron 3847, Intel(R) Core(TM) i5-4460 CPU@ 3.20GHz).

\subsection{An ODE problem with complicated coefficients}
{    
At first, let us consider the following ordinary differential equation (ODE) with complicated non-constant coefficients
}
\begin{equation}\label{Eq:eg2}
    x^2\frac{d^2u}{dx^2}+\sin\pi x\frac{du}{dx}+e^xu=\left(e^x+\pi\cos\pi x-\pi^2 x^2\right)\sin\pi x,
\end{equation}
subject to the boundary conditions
\begin{equation}\label{BC:eg2}
   u(0)=u(1)=0,
\end{equation}
which has the exact solution
\begin{equation}\label{Eq:eg2:solution}
   u(x)=\sin \pi x.
\end{equation}

According to \textbf{Theorem 2}, Eq. (\ref{Eq:eg2}) can be rewritten as
\begin{equation}\label{Eq:eg2:trans}
   \frac{d^2\left(b_2(x)u(x)\right)}{dx^2}+\frac{d\left(b_1(x)u(x)\right)}{dx}+b_0(x)u(x)=r(x),
\end{equation}
where
\begin{equation}
    \left\{
    \begin{aligned}
      b_0(x)&=e^x-\pi \cos\pi x+2, \\
      b_1(x)&=\sin \pi x-4x,  \\
     b_2(x)&=x^2, \\
     r(x)&=\left(e^x+\pi\cos\pi x-\pi^2 x^2\right)\sin\pi x.
    \end{aligned}
    \right.
\end{equation}

Considering the boundary condition $u(0)=0$ and the fact that $x=0$ is a first order zero point of $b_1(x)$ and a second order zero point of $b_2(x)$, we have
\begin{equation}
    \left\{
    \begin{aligned}
     &b_0(x)u(x)|_{x=0}=0,\\
     &b_1(x)u(x)|_{x=0}=0, \frac{d\left(b_1(x)u(x)\right)}{dx}\Big|_{x=0}=0, \\
     &b_2(x)u(x)|_{x=0}=0, \frac{d\left(b_2(x)u(x)\right)}{dx}\Big|_{x=0}=0, \frac{d^2\left(b_2(x)u(x)\right)}{dx^2}\Big|_{x=0}=0.
    \end{aligned}
    \right.
\end{equation}
So wavelet basis of $b_1(x)u(x)$ and $b_2(x)u(x)$ should be modified into $\tilde{\varphi}_{j,k}(x)$ and $\tilde{\tilde{\varphi}}_{j,k}(x)$, respectively, which are defined by Eqs. (\ref{Def:wavelet modification:tilde}) and (\ref{Def:wavelet modification:double-tilde}). Substituting these modified expressions into the wavelet approximation (\ref{Eq:ODE:wavelet}), we have
\begin{align}\label{Eq:eg2:wavelet}
   \sum_{k=0}^{2^j}&u\left(\frac{k}{2^j}\right)\left[b_2\left(\frac{k}{2^j}\right)\frac{d^2\tilde{\tilde{\varphi}}_{j,k}(x)}{dx^2}+b_1\left(\frac{k}{2^j}\right)\frac{d\tilde{\varphi}_{j,k}(x)}{dx}+b_0\left(\frac{k}{2^j}\right)  \varphi_{j,k}(x) \right]=0, \nonumber \\
\end{align}
where the boundary conditions can be embedded by directly substituting $u\left(0/2^j\right)=0$ and $u\left(2^j/2^j\right)=0$ into it.

Multiplying both sides of Eq. (\ref{Eq:eg2:wavelet}) by $\varphi_{j,l}(x)$ $(l=1, 2, ..., 2^j-1)$ respectively, and integrating over the interval $[0,1]$, we have a linear algebraic equation
\begin{equation}\label{Eq:eg2:wavelet-Galerkin}
   \left[\sum_{n=0}^2 \mathbf{A_n}^T \odot \mathbf{b_n}\right] \mathbf{u}=0.
\end{equation}
where
\begin{equation}
    \left\{
    \begin{aligned}
    &\mathbf{A_0}=\left\{a_{kl}=\Gamma_{k,l}^{j,0}=\int_{0}^{1}\varphi_{j,k}(\xi) \varphi_{j,l}(\xi)d\xi\right\}_{(2^j-1)\times (2^j+1)}  \\
    &\mathbf{A_1}=\left\{a_{kl}=\tilde{\Gamma}_{k,l}^{j,1}=\int_{0}^{1}\frac{d\tilde{\varphi}_{j,k}(\xi)}{d\xi} \varphi_{j,l}(\xi)d\xi\right\}_{(2^j-1)\times (2^j+1)}  \\
    &\mathbf{A_2}=\left\{a_{kl}=\tilde{\tilde{\Gamma}}_{k,l}^{j,2}=\int_{0}^{1}\frac{d^2\tilde{\tilde{\varphi}}_{j,k}(\xi)}{d\xi^2} \varphi_{j,l}(\xi)d\xi\right\}_{(2^j-1)\times (2^j+1)}  \\
    & \mathbf{b_{n}}=\left\{b_{k}=b_{n}\left(\frac{k}{2^j}\right)\right\}_{(2^j+1)\times 1}\\
    & \mathbf{u}=\left\{u_{k}=u\left(\frac{k}{2^j}\right)\right\}_{(2^j+1)\times 1}
    \end{aligned}
    \right.
\end{equation}
Separating the known parts that are relevant to $u(0/2^j)$ and $u(2^j/2^j)$ from the left side, it is very easy to solve the linear algebraic equation (\ref{Eq:eg2:wavelet-Galerkin}).

\begin{figure}[!h]
\centering
  \includegraphics[width=12cm]{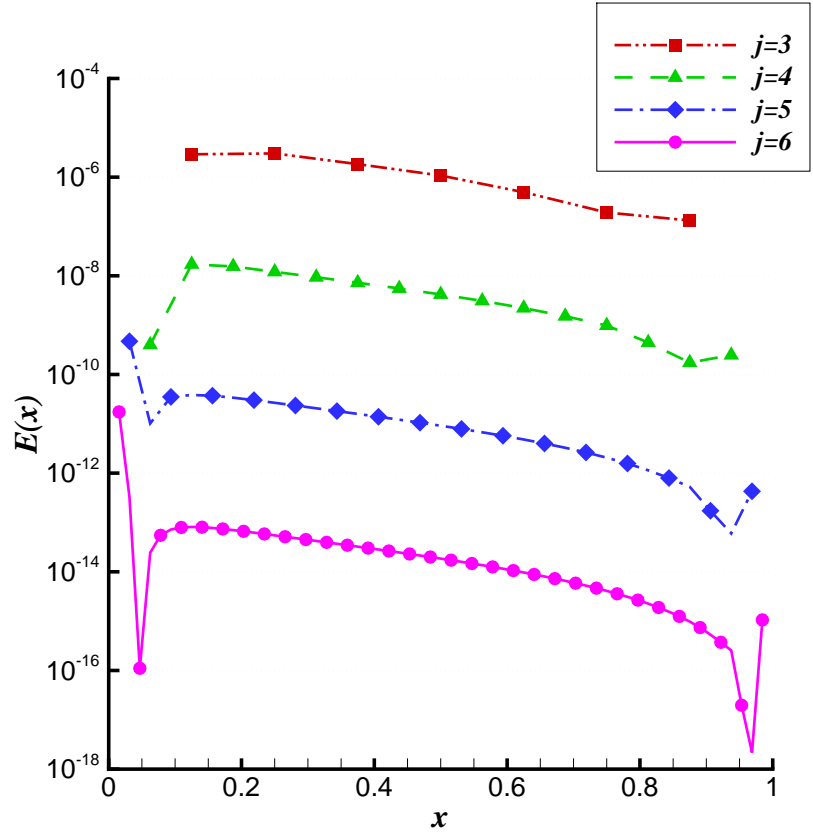}
      
  \renewcommand{\figurename}{Figure}
  \caption{The error distribution of the wavelet solution at different resolution levels for Eqs. (\ref{Eq:eg2}) and (\ref{BC:eg2}).
  Square symbol with double dot-dashed line:  $j=3$;
  delta symbol with dashed line: $j=4$;
  diamond symbol with dash-dotted line: $j=5$;
  circle symbol with solid line: $j=6$.
  }
\label{Fig:Eg2:ErrorDistribution}
\end{figure}

Define the averaged square error $ErrSQ$ and the error distribution $E(x)$ as follows:
\begin{equation}\label{Def:ErrSQ:1D}
    ErrSQ=\frac{1}{2^j+1} \sum_{k=0}^{2^j}\left[u_w(\frac{k}{2^j})-u_e(\frac{k}{2^j})\right]^2,
\end{equation}
\begin{equation}\label{Def:E(x)}
    E(x)=\big|u_w(x)-u_e(x)\big|^2,\ \ 0<x<1,
\end{equation}
where $u_e(x)$ denotes the exact solution and $u_w(x)$ denotes the wavelet solution.

\begin{table}[!h]
    
\tabcolsep 0pt \caption{The averaged square error and CPU time of wavelet solutions with different resolution levels of Eqs. (\ref{Eq:eg2}) and (\ref{BC:eg2}).} \label{Tab:Eg2:Accuracy and Efficiency} \vspace*{-12pt}
\begin{center}
\def\temptablewidth{1\textwidth}
{\rule{\temptablewidth}{1pt}}
\begin{tabular*}{\temptablewidth}{@{\extracolsep{\fill}}ccc}
$j$, resolution level  & $ErrSQ$ &  CPU time (sec.) \\
\hline
3	&	$1.1\times 10^{-6}$	&	0.39 \\
4	&	$4.7\times 10^{-9}$	&	0.41 \\
5	&	$2.6\times 10^{-11}$	&	0.55 \\
6	&	$8.0\times 10^{-14}$	&	1.01
\end{tabular*}
{\rule{\temptablewidth}{1pt}}
\end{center}
\end{table}

{    
The error distribution of the wavelet solution at different resolution levels are shown in Figure \ref{Fig:Eg2:ErrorDistribution}. The error distribution shows that the wavelet solutions agree very well with the exact solution. This verifies the validity of the generalized wavelet-Galerkin method. Also, by raising the resolution level $j$, the accuracy of the solution can be greatly improved. Table \ref{Tab:Eg2:Accuracy and Efficiency} clearly shows that by raising  the resolution level $j$ from 4 to 6, the averaged error quickly decreases from $10^{-6}$ to $10^{-14}$. More importantly, the advantages of this method are revealed. Note that the connection coefficients $\Gamma_{k,l}^{j,n}$ are independent of the non-constant coefficients $x^2$, $\sin \pi x$ and $e^x$ in Eq. (\ref{Eq:eg2}), in other words, the connection coefficients are independent of the differential equations at all.  Therefore,  it is very easy to build a database for data sharing, which of course makes this method very efficient and adaptive. For example, in case of the resolution level $j=4$, only 0.39 seconds CPU time is used to gain the wavelet solution. Even in case of the resolution level $j=6$, it just takes 1.01 seconds CPU time.  So, this approach is computationally rather efficient.
}

\subsection{An ODE problem with special functions as non-constant coefficients}

Let us further consider a BVP which contains non-constant coefficients expressed by special functions
\begin{equation}\label{Eq:eg3}
   \Gamma(x+1)\frac{d^2u(x)}{dx^2}+\left[\Gamma(x+1)+\Gamma'(x+1)\right]\frac{du(x)}{dx}=-\Gamma'(x+1)-\Gamma''(x+1),
\end{equation}
subject to the boundary conditions
\begin{equation}\label{BC:eg3}
   u(0)=0,\ u(1)=0,
\end{equation}
which has the exact solution
\begin{equation}\label{Eq:eg3:solution}
   u(x)=\ln \frac{1}{\Gamma(x+1)},
\end{equation}
where
\begin{equation}\label{Eq:GammaFunction}
   \Gamma(x)=\int_0^{+\infty}\xi^{x-1}e^{-\xi}d\xi.
\end{equation}

According to \textbf{Theorem 2}, Eq. (\ref{Eq:eg3}) can be rewritten as
\begin{equation}\label{Eq:eg:trans}
   \frac{d^2\left(b_2(x)u(x)\right)}{dx^2}+\frac{d\left(b_1(x)u(x)\right)}{dx}+b_0(x)u(x)=r(x),
\end{equation}
where
\begin{equation}
    \left\{
    \begin{aligned}
     b_0(x)&=-\Gamma'(x+1),  \\
     b_1(x)&=\Gamma(x+1)-\Gamma'(x+1),  \\
     b_2(x)&=\Gamma(x+1),  \\
     r(x)&=-\Gamma'(x+1)-\Gamma''(x+1).
    \end{aligned}
    \right.
\end{equation}

\begin{figure}[!h]
\centering
  \includegraphics[width=12cm]{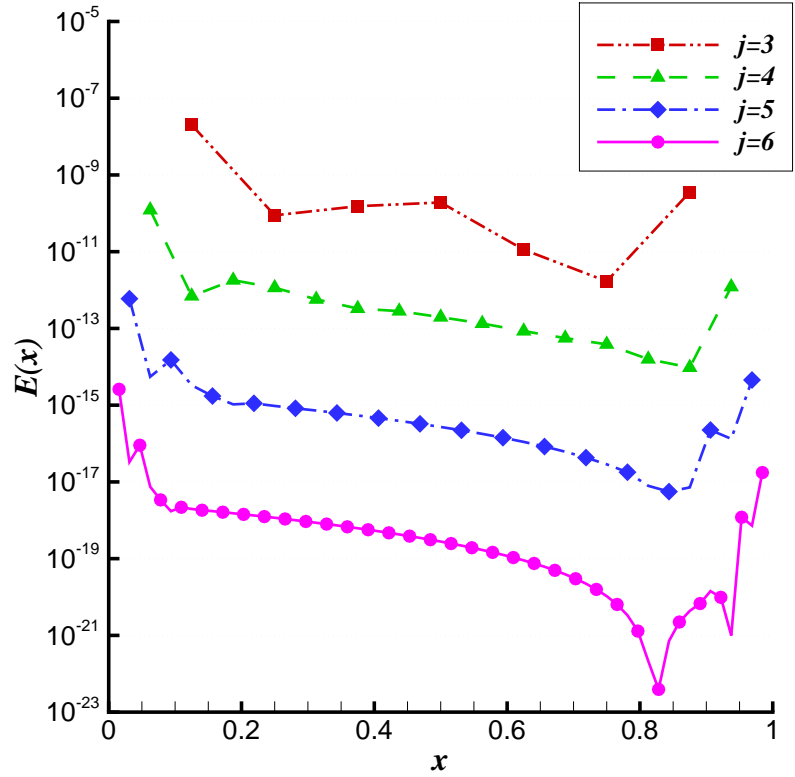}
      
  \renewcommand{\figurename}{Figure}
  \caption{The error distribution of the wavelet solution at different resolution levels for Eqs. (\ref{Eq:eg3}) and (\ref{BC:eg3}).
  Square symbol with double dot-dashed line:  $j=3$;
  delta symbol with dashed line: $j=4$;
  diamond symbol with dash-dotted line: $j=5$;
  circle symbol with solid line: $j=6$.
  }
\label{Fig:Eg3:ErrorDistribution}
\end{figure}

For this problem, as $x=0$ and $x=1$ are not zero points of $b_0(x)$, $b_1(x)$ and $b_2(x)$, no modification on the wavelet basis is necessary, so the boundary conditions (\ref{BC:eg3}) can be embedded into the wavelet expressions of $b_0(x) u(x)$, $b_1(x) u(x)$ and $b_2(x) u(x)$ by directly substituting $u(0)=u(1)=0$ into them. Thus, the solution of this boundary value problem can be obtained by substituting these coefficient functions into the generalized wavelet-Galerkin equation (\ref{Eq:ODE:wavelet-Galerkin}) and solving it.

{    
The error distribution of the wavelet solutions at different resolution levels  are presented in Figure \ref{Fig:Eg3:ErrorDistribution}. Obviously, the error is rather small, which shows that the wavelet solutions agree very well with the exact solution. This once again verifies the validity of the generalized wavelet-Galerkin method.
}

\begin{table}[!h]
\tabcolsep 0pt \caption{The averaged square error and CPU time of wavelet solutions at different resolution levels of Eqs. (\ref{Eq:eg3}) and (\ref{BC:eg3}).} \label{Tab:Eg3:Accuracy and Efficiency} \vspace*{-12pt}
\begin{center}
\def\temptablewidth{1\textwidth}
{\rule{\temptablewidth}{1pt}}
\begin{tabular*}{\temptablewidth}{@{\extracolsep{\fill}}ccc}
$j$, resolution level  & $ErrSQ$ &  CPU time (sec.) \\
\hline
3	&	$2.3\times 10^{-09}$	&	0.31 \\
4	&	$7.6\times 10^{-12}$	&	0.37 \\
5	&	$1.9\times 10^{-14}$	&	0.54 \\
6	&	$4.3\times 10^{-17}$	&	0.95
\end{tabular*}
{\rule{\temptablewidth}{1pt}}
\end{center}
\end{table}

Table \ref{Tab:Eg3:Accuracy and Efficiency} presents the averaged square error and the corresponding used CPU time for the wavelet solutions at different resolution levels. Obviously, the generalized wavelet-Galerkin method possesses very high efficiency. With the resolution level $j$ varying from 3 to 6, all of the wavelet solutions can be obtained in less than one second while the computational accuracy rises significantly. At the resolution level $j=6$, the averaged square error is rather small, say, $10^{-17}$. Moreover, it can be seen in Figure \ref{Fig:Eg3:ErrorDistribution} that almost all the averaged square errors on the whole interval is smaller than $10^{-15}$ at the resolution level $j=6$. It means that the error is rather small on the whole interval. These results not only verify the validity of the generalized wavelet-Galerkin method, but also illustrate the good convergence property, high accuracy and high efficiency of computation as well as the great potential of this method.

\subsection{A PDE problem}

{    Furthermore}, let us consider a two-dimensional boundary value problem governed by
\begin{equation}\label{Eq:eg4}
    \sqrt{x^2+y^2+1}\left(\frac{\partial^2u}{\partial x^2}+\frac{\partial^2u}{\partial x^2}\right)+\frac{1}{ \sqrt{x^2+y^2+1}}\left(x\frac{\partial u}{\partial x}+y\frac{\partial u}{\partial y}-2u\right)=0,
\end{equation}
subject to the boundary conditions
\begin{equation}\label{BC:eg4}
    \left\{
    \begin{aligned}
     u(0,y)=\sqrt{1+y^2},\ u(1,y)=\sqrt{2+y^2},\\
    u(x,0)=\sqrt{1+x^2},\ u(x,1)=\sqrt{2+x^2},
	\end{aligned}
    \right.
\end{equation}
which has the exact solution
\begin{equation}\label{Eq:eg4:solution}
   u(x,y)=\sqrt{x^2+y^2+1}.
\end{equation}

According to \textbf{Theorem 4}, Eq. (\ref{Eq:eg4}) can be rewritten as
\begin{equation}\label{Eq:eg4:trans}
   \frac{\partial^2\left(b_{20}u\right)}{\partial x^2}+\frac{\partial^2\left(b_{02}u\right)}{\partial y^2}+\frac{\partial^2\left(b_{11}u\right)}{\partial x \partial y}+\frac{\partial\left(b_{10}u\right)}{\partial x}+\frac{\partial\left(b_{01}u\right)}{\partial y}+b_{00}u=0,
\end{equation}
where
\begin{equation}
    \left\{
    \begin{aligned}
     b_{20}(x,y)&=\sqrt{x^2+y^2+1}, \\
     b_{02}(x,y)&=\sqrt{x^2+y^2+1},  \\
     b_{11}(x,y)&=0,  \\
     b_{10}(x,y)&=-\frac{x}{\sqrt{x^2+y^2+1}},  \\
     b_{01}(x,y)&=-\frac{y}{\sqrt{x^2+y^2+1}},  \\
     b_{00}(x,y)&=-\frac{2}{\sqrt{x^2+y^2+1}}
	\end{aligned}
    \right.
\end{equation}

\begin{figure}[!h]
\centering
  \includegraphics[width=12cm]{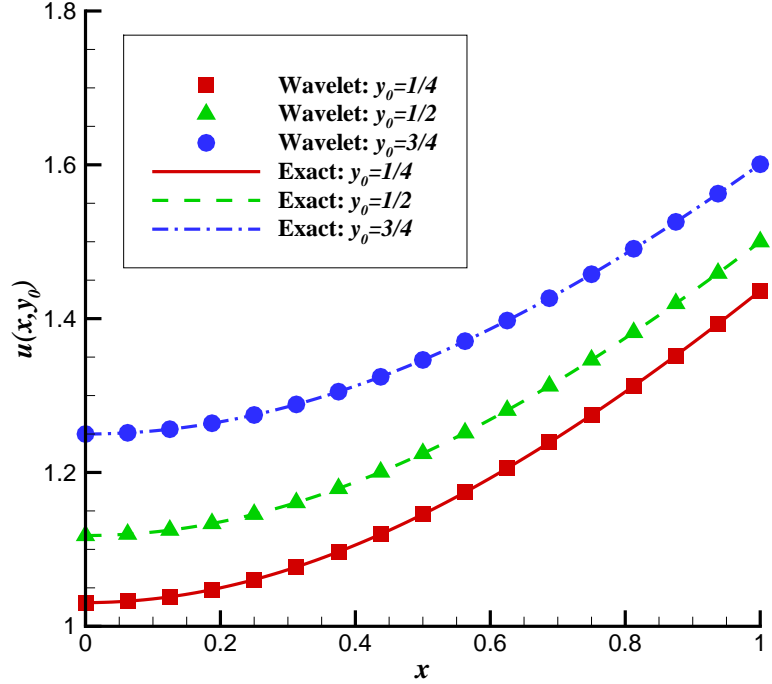}
  \renewcommand{\figurename}{Figure}
  \caption{Comparison of section curves of the wavelet solution at the resolution level $j=4$ with the exact solution of Eqs. (\ref{Eq:eg4}) and (\ref{BC:eg4}).
  Lines: exact solutions;
  square symbol: wavelet solution at $y_0=1/4$;
  delta symbol: wavelet solution at $y_0=1/2$;
  circle symbol: wavelet solution at $y_0=3/4$;
  }
\label{Fig:Eg4:SolutionComparison_Sections}
\end{figure}

\begin{figure}[!h]
\centering
  \includegraphics[width=12cm]{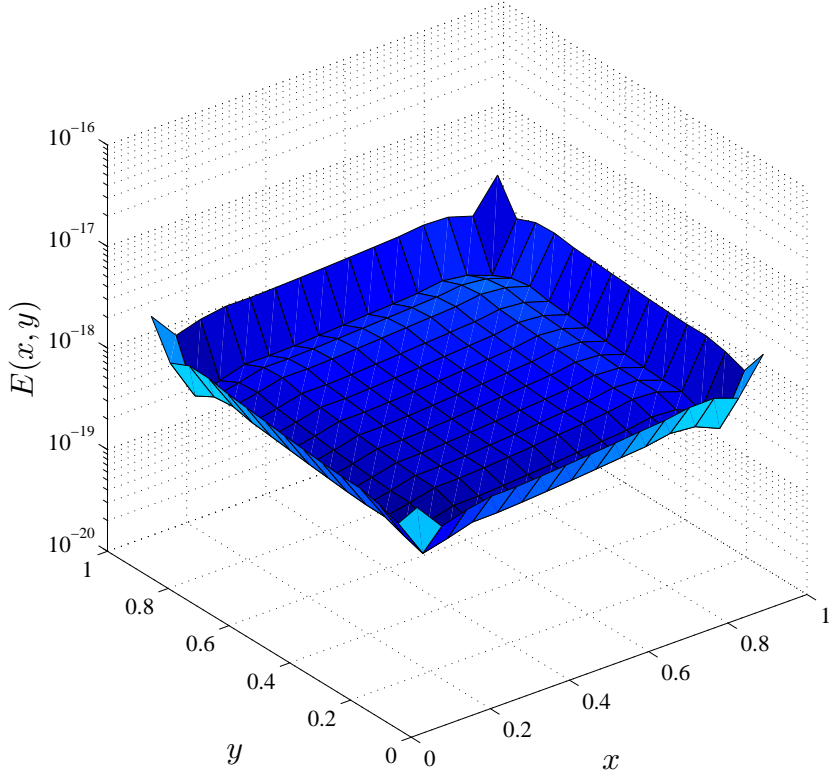}
      
  \renewcommand{\figurename}{Figure}
  \caption{The error distribution of the wavelet solution at the resolution level $j=4$ for Eqs. (\ref{Eq:eg4}) and (\ref{BC:eg4}).
  }
\label{Fig:Eg4:ErrorDistribution}
\end{figure}

For this problem, Eq. (\ref{Eq:eg4:trans}) can be approximated by
\begin{align}\label{Eq:eg4:wavelet}
   &\sum_{k=0}^{2^j}\sum_{l=0}^{2^j}u\left(\frac{k}{2^j},\frac{l}{2^j}\right)\left[b_{20}\left(\frac{k}{2^j},\frac{l}{2^j}\right)\frac{d^2\varphi_{j,k}(x)}{dx^2}\varphi_{j,l}(y)+b_{02}\left(\frac{k}{2^j},\frac{l}{2^j}\right)\varphi_{j,k}(x)\frac{d^2\varphi_{j,l}(y)}{dy^2}\right] \nonumber \\
  &+\sum_{k=0}^{2^j}\sum_{l=0}^{2^j}u\left(\frac{k}{2^j},\frac{l}{2^j}\right)\left[b_{10}\left(\frac{k}{2^j},\frac{l}{2^j}\right)\frac{d\varphi_{j,k}(y)}{dx}\varphi_{j,l}(y)+b_{01}\left(\frac{k}{2^j},\frac{l}{2^j}\right)\varphi_{j,k}(x)\frac{d\varphi_{j,l}(y)}{dy}\right] \nonumber \\
   &+\sum_{k=0}^{2^j}\sum_{l=0}^{2^j}u\left(\frac{k}{2^j},\frac{l}{2^j}\right)b_{00}\left(\frac{k}{2^j},\frac{l}{2^j}\right)\varphi_{j,k}(x)\varphi_{j,l}(y)=0,
\end{align}
The boundary conditions (\ref{BC:eg4}) can be embedded by directly substituting
\begin{equation}\label{BCdiscrete:eg4}
    \left\{
    \begin{aligned}
     u\left(\frac{k}{2^j},0\right)=\sqrt{1+\left(\frac{k}{2^j}\right)^2}, \ u\left(\frac{k}{2^j},1\right)=\sqrt{2+\left(\frac{k}{2^j}\right)^2}, \\
     u\left(0,\frac{l}{2^j}\right)=\sqrt{1+\left(\frac{l}{2^j}\right)^2},\ u\left(1,\frac{l}{2^j}\right)=\sqrt{2+\left(\frac{l}{2^j}\right)^2}
	\end{aligned}
    \right.
\end{equation}
into wavelet approximation (\ref{Eq:eg4:wavelet}) for $k,l=0,1,2,...,2^j$.

 Thus, we have
\begin{equation}\label{Eq:eg4:wavelet-Galerkin}
   \left[\sum_{m=0}^2\sum_{n=0}^{2-m}\left(\mathbf{A_m}\otimes \mathbf{A_n}\right)^T \odot rvec\left(\mathbf{B_{mn}}\right)\right] rvec\left(\mathbf{U}\right)=\mathbf{0},
\end{equation}
where the matrixes are defined by Eq. (\ref{Def:Matrixes:2D}).

Separating the known parts relevant to the boundary values from the left side and solving this liner algebraic equation, the solution of the two-dimensional BVP is obtained. The three section curves at $y=1/4$, $y=1/2$ and $y=3/4$ at the resolution level $j=4$ are shown in Figure \ref{Fig:Eg4:SolutionComparison_Sections}. The corresponding averaged squared error is only $7.5\times 10^{-19}$ as presented in Table \ref{Tab:Eg4:Accuracy and Efficiency}. {    Also, Figure \ref{Fig:Eg4:ErrorDistribution} shows that the error is evenly distributed.}Obviously, the wavelet approximation agrees very well with the exact solution, which further indicates the validity of the generalized wavelet-Galerkin method in two-dimensional boundary value problems.

\begin{table}[!h]
\tabcolsep 0pt \caption{The averaged square error and CPU time of wavelet solutions with different resolution levels of Eqs. (\ref{Eq:eg4}) and (\ref{BC:eg4}).} \label{Tab:Eg4:Accuracy and Efficiency} \vspace*{-12pt}
\begin{center}
\def\temptablewidth{1\textwidth}
{\rule{\temptablewidth}{1pt}}
\begin{tabular*}{\temptablewidth}{@{\extracolsep{\fill}}ccc}
$j$, resolution level  & $ErrSQ$ &  CPU time (sec.) \\
\hline
3	&	$8.7\times 10^{-19}$	&	0.36 \\
4	&	$7.5\times 10^{-19}$	&	0.44 \\
5	&	$6.9\times 10^{-19}$	&	0.74 \\
6	&	$6.6\times 10^{-19}$	&	3.53
\end{tabular*}
{\rule{\temptablewidth}{1pt}}
\end{center}
\end{table}

{    
\section{Application in nonlinear problems}
For simplicity, the basic idea of the generalized wavelet-Galerkin method is mainly described by linear cases. However, the significance of the  generalized wavelet-Galerkin method extends far beyond solving linear differential equations with non-constant coefficients. Actually, the generalized wavelet-Galerkin method is also a powerful tool to solve complicated nonlinear problems with non-constant coefficients.

To illustrate it more clearly, let us consider a nonlinear BVP which contains non-constant coefficients as follows:
\begin{equation}\label{Eq:eg5}
   \frac{d^2u}{dx^2} - 2(x+1) e^{2 x}u\frac{du}{dx}+ 2\ln u+\left(\sin\pi x + 2 e^x\right)\frac{du}{dx} + \left(\sin\pi x -1\right)u =0,
\end{equation}
subject to the boundary conditions
\begin{equation}\label{BC:eg5}
   u(0)=1,\ u(1)=\frac{1}{e},
\end{equation}
in which $- 2(x+1) e^{2 x}u\frac{du}{dx}$ is a nonlinear term with non-constant coefficient and $2\ln u$ is a transcendental nonlinear term.

It is known that one exact solution of the nonlinear BVP is
\begin{equation}\label{Eq:eg5:solution}
   u(x)= e^{-x}.
\end{equation}
However, as it is a complicated nonlinear boundary value problem, there might be more than one solution.

Similar to the basic ideas of \textbf{Theorem 1} and \textbf{Theorem 2}, Eq. (\ref{Eq:eg5}) can be rewritten as
\begin{equation}\label{Eq:eg5:trans}
   \frac{d^2\left(b_{2}(x)u\right)}{dx^2}+ \frac{d\left(b_{1}(x)u\right)}{dx}+b_{0}(x)u+\frac{d\left(c_{1}(x)u^2\right)}{dx}+c_{0}(x)u^2 + d_{0}(x)\ln u=0,
\end{equation}
where
\begin{equation}
    \left\{
    \begin{aligned}
     b_{0}(x)&=\sin(\pi x) - \pi \cos(\pi x) - 2 e^x -1 ,  \\
     b_{1}(x)& = \sin(\pi x) + 2 e^x,  \\
     b_{2}(x)& = 1,  \\
     c_{0}(x)& = (2x+3)e^{2x},  \\
     c_{1}(x)& = -(x+1)e^{2x},\\
     d_{0}(x)& = 2.
    \end{aligned}
    \right.
\end{equation}

Using the generalized wavelet-Galerkin method, we have:
\begin{equation}\label{Eq:eg5:wavelet-Galerkin}
   \left[\sum_{n=0}^2 \mathbf{A_n}^T \odot \mathbf{b_{n}}\right] \mathbf{u} + \left[\sum_{n=0}^1 \mathbf{A_n}^T \odot \mathbf{c_{n}}\right] \mathbf{u} \circ \mathbf{u}+\left[\mathbf{A_0}^T \odot \mathbf{d_{0}}\right] \ln( \mathbf{u})=0.
\end{equation}
Separating the known parts relevant to the boundary values from the left side and solving this nonlinear algebraic equation, the solutions of the nonlinear BVP can be conveniently obtained.

What's of special interest is that by the generalized wavelet-Galerkin method, more than one solution is obtained, as is shown in Figure \ref{Fig:Eg5:SolutionComparison}. Obviously, wavelet solution \uppercase\expandafter{\romannumeral1} (the lower branch) agrees very well with the exact solution, and the wavelet solution \uppercase\expandafter{\romannumeral2} (the upper branch) is undoubtedly a newly found solution, though the corresponding exact one of which is unknown. It verifies that the generalized wavelet-Galerkin method is a powerful tool to solve nonlinear differential equations with non-constant coefficients. Especially, it is possible to find new solutions of multi-solution problems. It means that the generalized wavelet-Galerkin method possesses great potential for application in nonlinear problems.

\begin{figure}[!h]
\centering
  \includegraphics[width=12cm]{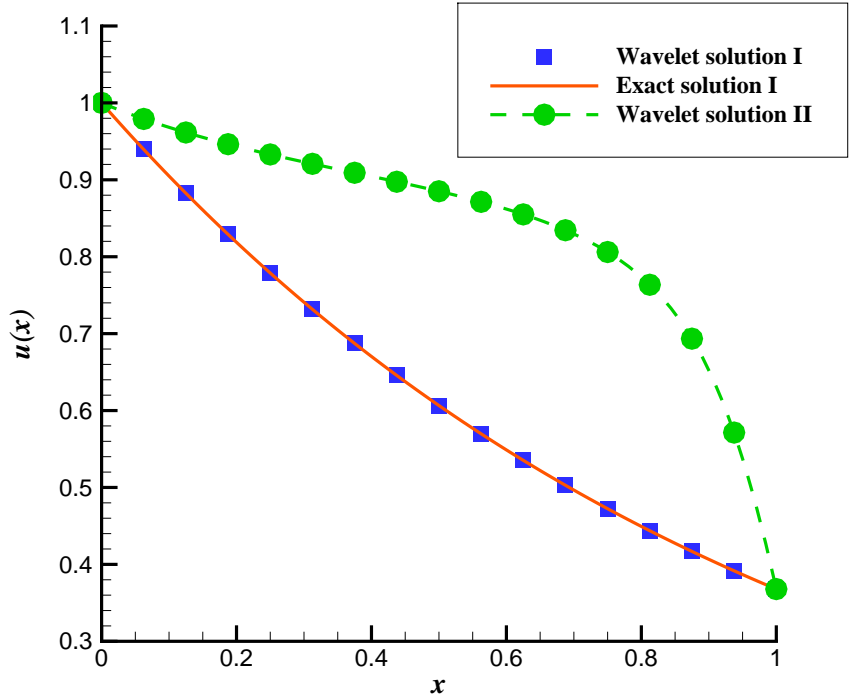}
      
  \renewcommand{\figurename}{Figure}
  \caption{Comparison of the wavelet solution at the resolution level $j=4$ with the exact solution, and the newly found wavelet solution for Eqs. (\ref{Eq:eg5}) and (\ref{BC:eg5}).
  Solid line: the known exact solution;
  square symbol: wavelet solution I;
  circle symbol with dashed line: wavelet solution II (newly found).
  }
\label{Fig:Eg5:SolutionComparison}
\end{figure}

\begin{figure}[!h]
\centering
  \includegraphics[width=12cm]{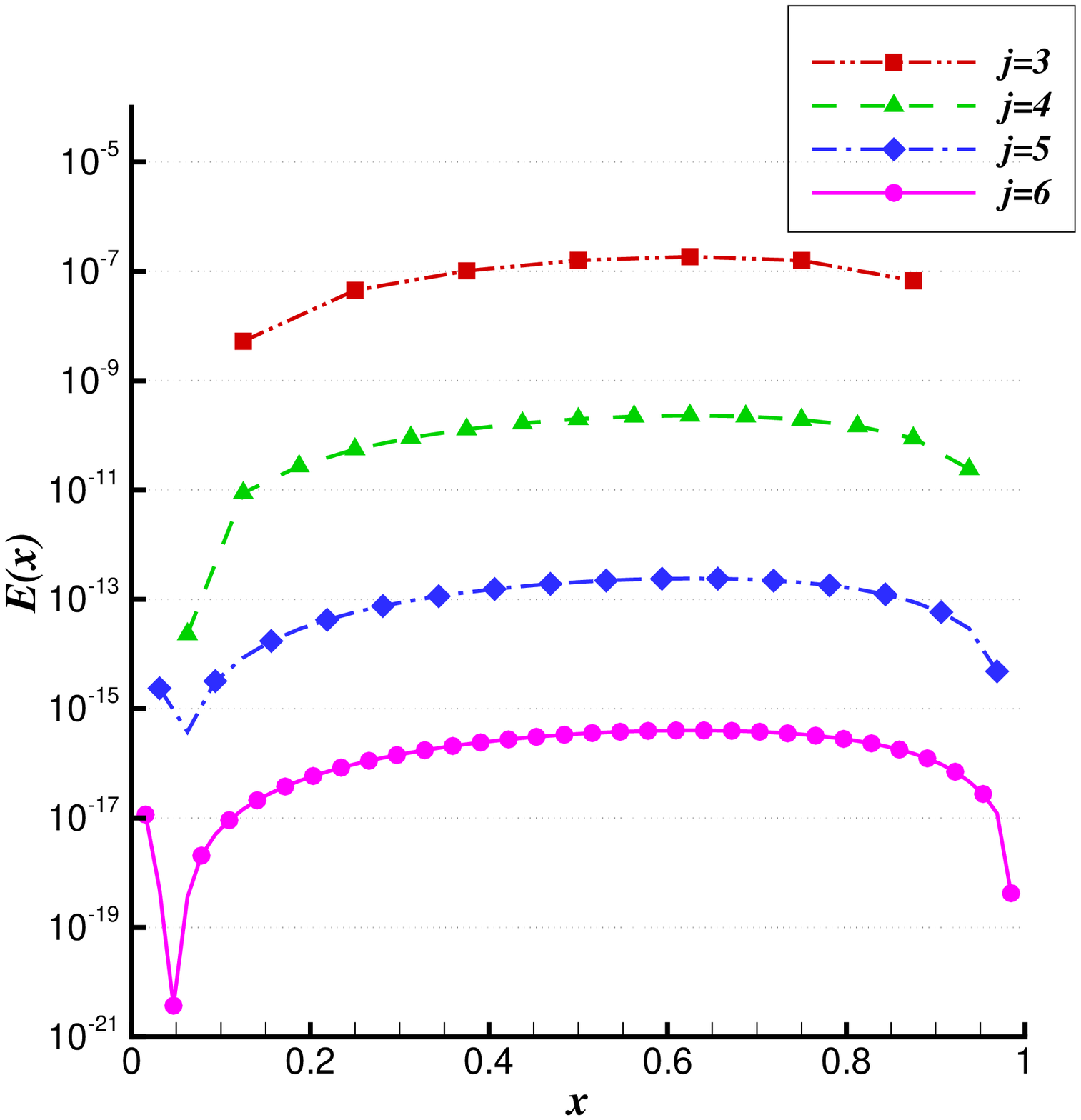}
      
  \renewcommand{\figurename}{Figure}
  \caption{The error distribution of the wavelet solution I (the lower branch) at different resolution levels for Eqs. (\ref{Eq:eg5}) and (\ref{BC:eg5}).
  Square symbol with double dot-dashed line:  $j=3$;
  delta symbol with dashed line: $j=4$;
  diamond symbol with dash-dotted line: $j=5$;
  circle symbol with solid line: $j=6$.
  }
\label{Fig:Eg5:ErrorDistribution}
\end{figure}

\begin{table}[!h]
    
\tabcolsep 0pt \caption{The averaged square error and CPU time of wavelet solutions with different resolution levels of Eqs. (\ref{Eq:eg5}) and (\ref{BC:eg5}).} \label{Tab:Eg5:Accuracy and Efficiency} \vspace*{-12pt}
\begin{center}
\def\temptablewidth{1\textwidth}
{\rule{\temptablewidth}{1pt}}
\begin{tabular*}{\temptablewidth}{@{\extracolsep{\fill}}ccc}
$j$, resolution level  & $ErrSQ$ &  CPU time (sec.) \\
\hline
3	&	$8.0\times 10^{-8}$	&	0.72 \\
4	&	$1.1\times 10^{-10}$	&	0.80 \\
5	&	$1.2\times 10^{-13}$	&	0.95 \\
6	&	$1.9\times 10^{-16}$	&	1.52
\end{tabular*}
{\rule{\temptablewidth}{1pt}}
\end{center}
\end{table}

To verify the accuracy of the generalized wavelet-Galerkin method for nonlinear cases, the error distribution curves of the wavelet solution I (with known exact solution) at different resolution levels are shown in Figure \ref{Fig:Eg5:ErrorDistribution}. It strongly verifies the validity as well as the high accuracy of the generalized wavelet-Galerkin method. In addition, the averaged errors at different resolution levels and the corresponding CPU time are present in Table \ref{Tab:Eg5:Accuracy and Efficiency}. It shows that with the resolution level raising from 3 to 6, the averaged error decreases from $10 ^ {-8}$ to $10 ^ {-16}$, while the CPU time is less than 2 seconds, which strongly indicates that the generalized wavelet-Galerkin method possesses both high accuracy and high efficiency not only for linear problems, but also for nonlinear ones. This once again shows the great potential of the generalized wavelet-Galerkin method.}

\section{Concluding remarks and discussions}

The computation of connection coefficients plays an important role in the wavelet-Galerkin method. However, differential equations with non-constant coefficients lead to very complicated connection coefficients, which seriously limits the application of the traditional wavelet-Galerkin method based on compactly supported wavelets. In this paper, a generalized wavelet-Galerkin method is proposed to overcome these  restrictions of the traditional one.

Some basic theorems are proposed and proved. Based on these theorems, the generalized wavelet-Galerkin method is developed. In the frame of the generalized wavelet-Galerkin method, complicated connection coefficients are always independent of the type of linear differential equations. Therefore, it is easy to build a database of the basic connection coefficients for data sharing, because no new connection coefficient arises no matter what kinds of non-constant coefficients are contained in a differential equation. That makes the generalized wavelet-Galerkin method very efficient and adaptive for various types of boundary value problems. The generalized wavelet-Galerkin equations are derived in a very concise form, which are of great significance both theoretically and practically. Using the generalized Coiflet-type wavelet, several linear boundary value problems governed by ordinary and partial differential equations with non-constant coefficients are solved by the generalized wavelet-Galerkin method. All of these examples verify the validity as well as the advantages of the generalized wavelet-Galerkin method. {    More importantly, the nonlinear application indicates that the generalized wavelet-Galerkin method not only keeps the same advantages for nonlinear problems, but also possesses the ability to find new solutions for multi-solution problems.}

Note that the generalized wavelet-Galerkin method possesses very great potential for extension. Although only one-dimensional and two-dimensional boundary value problems are illustrated  in this paper, the generalized wavelet-Galerkin method is valid for higher dimensional problems. Actually, the generalized wavelet-Galerkin equation (\ref{Eq:PDE:wavelet-Galerkin}) can be extended to $K$-dimensional BVPs:
\begin{align}\label{Eq:high-dimensional:wavelet-Galerkin}
  & \left[\sum_{n_1=0}^{N_1}\sum_{n_2=0}^{N_2}...\sum_{n_K=0}^{N_K}\left(\mathop{\bigotimes} \limits_{k=1}^K\mathbf{A_{n_k}} \right)^T \odot rvec\left(\mathbf{B_{n_1n_2...n_K}}\right)\right] rvec\left(\mathbf{U}\right) \nonumber \\
   &\quad \quad  \quad  \quad  \quad  \quad  \quad  \quad  \quad \quad  \quad  \quad  \quad  \quad \quad \quad \quad =\left(\mathop{\bigotimes}_{k=1}^K \mathbf{A_0}\right)^T rvec\left(\mathbf{R}\right),
\end{align}
where $\mathop{\bigotimes} \limits_{k=1}^K\mathbf{A_{n_k}}=\mathbf{A_{n_1}}\otimes\mathbf{A_{n_2}}\otimes...\otimes\mathbf{A_{n_K}}$, and the operator $rvec\left(\cdot\right)$ pulls a $K$-dimensional matrix into a vector dimension by dimension. {    Of course, this equation can be further extended and applied to nonlinear cases.} On the other hand, though only the generalized Coiflet-type wavelet is discussed in this paper, theoretically it is possible to extend the generalized wavelet-Galerkin method for other kinds of wavelets with interpolation or quasi-interpolation properties.

Note that Yang and Liao \cite{zcYang2017-CNSNS-A, zcYang2017-CNSNS-B} successfully combined the wavelet with the homotopy analysis method (HAM) \cite{liao2004-AMC, Liaobook, liaobook2, KV2012}, an analytic approximation method for highly nonlinear problems developed by Liao \cite{LiaoPhd} in 1992.  Unlike perturbation methods, the HAM has nothing to do with the small/large physical parameters.  Besides, different from all other approximation techniques, the HAM provides us a convenient way to guarantee the convergence of solution series.   As illustrated by Zhong and Liao \cite{Zhong2017-SCPMA},  perturbation methods are often only special cases of the HAM, but the HAM still works well even if perturbation methods fail.  More importantly, as a new analytic method, the HAM can provide us something completely new/different \cite{xu2012JFM, BookChap3-Liao2015, Liao2016JFM}.   So, we can combine the generalized Wavelet-Galerkin method with the HAM to solve highly nonlinear equations.    

\section*{Acknowledgment}

This work is partially supported by  the National Natural Science Foundation of China  (Approval No. 11272209 and 11432009).

\bibliographystyle{model3-num-names}
\bibliography{RefGWGM}

\end{document}